\documentclass[a4paper]{amsart}
\usepackage[french,english]{babel}
\usepackage{amssymb} 
\usepackage[T1]{fontenc}
\usepackage[latin1]{inputenc}
\usepackage{amsfonts}
\usepackage{amsxtra}
\usepackage{ae}
\usepackage{slashed}
\usepackage{pdfsync}
\usepackage[all]{xy}
\usepackage{enumerate}
\usepackage{verbatim}
\usepackage{pgf,tikz}
\usepackage{url}
\include{diagram}
\usetikzlibrary{arrows,patterns}
\usepackage{color}
\usepackage[linkbordercolor={0.8 0.9 0.9}, citebordercolor={0.85 0.85 0.85}]{hyperref}

\newcommand*{\ket}{\rangle}
\newcommand*{\bra}{\langle}
\newcommand*{\ad}{\mathsf{ad}}

\newcommand*{\M}{\mathcal{M}}
\newcommand*{\N}{\mathcal{N}}
\newcommand*{\E}{\mathcal{E}}
\newcommand*{\F}{\mathcal{F}}

\renewcommand*{\H}{\mathcal{H}}

\renewcommand*{\O}{\mathcal{O}}

\newcommand*{\CF}{\mathfrak{C}}
\newcommand*{\DF}{\mathfrak{D}}
\newcommand*{\res}{\mathsf{res}}
\newcommand*{\ext}{\mathsf{ext}}

\newcommand*{\Poly}{\mathcal{O}}
\renewcommand*{\max}{\mathsf{f}}
\newcommand*{\red}{\mathsf{r}}
\newcommand*{\Irr}{\mathsf{Irr}}

\newcommand*{\hit}{\rightharpoonup}

\newcommand*{\KH}{\mathbb{K}}
\newcommand*{\LH}{\mathbb{L}}

\newcommand{\roots}{\mathbf{Q}}
\newcommand{\simpleroots}{\Sigma}
\newcommand{\weights}{\mathbf{P}}

\DeclareMathOperator{\End}{End}
\DeclareMathOperator{\tr}{tr}

\DeclareMathOperator{\Hom}{Hom}

\DeclareMathOperator{\im}{im}

\DeclareMathOperator{\id}{id}

\DeclareMathOperator{\HC}{HC}

\newcommand{\SL}{\operatorname{SL}}
\newcommand{\SU}{\operatorname{SU}}

\numberwithin{equation}{section}
\theoremstyle{change}
\newtheorem{theorem}{Theorem}[section]
\newtheorem{prop}[theorem]{Proposition}
\newtheorem{lemma}[theorem]{Lemma}

\newtheorem{definition}[theorem]{Definition}

\begin{document}

\title[Plancherel formula for complex quantum groups]{The Plancherel formula for complex semisimple quantum groups \\[2ex]
La formule de Plancherel pour les groupes quantiques semi-simples complexes}

\author{Christian Voigt}
\address{School of Mathematics and Statistics \\
         University of Glasgow \\
         University Place \\
         Glasgow G12 8QQ \\
         United Kingdom 
}

\email{christian.voigt@glasgow.ac.uk}

\author{Robert Yuncken}
\address{Laboratoire de Math\'ematiques \\ 
         Universit\'e Blaise Pascal \\
         Complexe universitaire des C\'ezeaux \\
         63177 Aubi\`ere Cedex \\
         France 
}

\email{Robert.Yuncken@math.univ-bpclermont.fr}


\thanks{The first author was supported by the Polish National Science Centre grant no. 2012/06/M/ST1/00169. 
The second author was supported by the project SINGSTAR of the Agence Nationale de la Recherche, ANR-14-CE25-0012-01, and by the CNRS PICS project OpPsi.}


\maketitle

\selectlanguage{english} 
\begin{abstract}
We calculate the Plancherel formula for complex semisimple quantum groups, that is, Drinfeld doubles of $ q $-deformations of 
compact semisimple Lie groups. As a consequence we obtain a concrete description of their associated reduced group $ C^* $-algebras. 
The main ingredients in our proof are the Bernstein-Gelfand-Gelfand complex and the Hopf trace formula.  
\end{abstract}

\selectlanguage{french} 
\begin{abstract}
Nous calculons la formule de Plancherel pour les groupes quantiques semi-simples complexes, c'est-\`a-dire les doubles de Drinfeld des $q$-d\'eformations des groupes de Lie semi-simples compacts.  En cons\'equence nous obtenons une description concr\`ete des $C^*$-alg\`ebres r\'eduites de groupe associ\'ees.  Les ingr\'edients principaux dans notre preuve sont le complexe de Bernstein-Gelfand-Gelfand et la formule des traces de Hopf.
\end{abstract}

\selectlanguage{english}

\section{Introduction}

Complex semisimple quantum groups are locally compact quantum groups which were constructed and first studied by Podle\'s and Woronowicz \cite{PWlorentz}. 
They are defined as Drinfeld doubles of $ q $-deformations of compact semisimple Lie groups, and can be viewed as deformations of the corresponding complex 
Lie groups in a natural way. Motivated by physical considerations, Podle\'s and Woronowicz focussed mainly on the case of the quantum Lorentz group, 
that is, the Drinfeld double of $ \SU_q(2) $. It became clear later that the theory of more general complex semisimple quantum groups is linked with a range 
of seemingly unrelated problems in noncommutative geometry, operator $ K $-theory, and the theory of $ C^* $-tensor categories and subfactors, see for 
instance \cite{Aranocomparison}, \cite{NVpoincare}, \cite{VYfredholm}. 

In the present paper we study the reduced unitary dual of complex semisimple quantum groups, and our main result is an explicit computation 
of the Plancherel formula. This generalizes work of Buffenoir and Roche \cite{BRLorentz} on the quantum Lorentz group. 
The formula we obtain can be interpreted as a deformation of the Plancherel formula for the corresponding classical groups, 
but our method of proof is completely different. 

We note that the abstract Plancherel theorem for locally compact quantum groups 
was established by Desmedt \cite{Desmedtthesis}, in analogy 
to the classical theory. In the case of complex semisimple quantum groups the Plancherel theorem involves so-called 
Duflo-Moore operators because the dual Haar weights fail to be traces. This is analogous to the situation for non-unimodular locally compact groups 
treated by Duflo and Moore in \cite{DufloMoorenonunimodular}. A classical locally compact group is unimodular if and only if the Haar weight on 
its group $ C^* $-algebra is tracial. In the quantum setting, traciality of the dual Haar weights implies unimodularity, 
but the converse does not hold in general. This is a well-known phenomenon which already shows up in the theory of compact quantum groups. 

Given a locally compact group or quantum group,  
a key problem is to calculate the Plancherel formula, that is, to determine 
explicitly the Plancherel measure and Duflo-Moore operators in terms of a given parametrization of the unitary dual. 
Before describing our proof strategy in the case of complex quantum groups, let us briefly recall the approach to compute the Plancherel formula for 
classical complex semisimple Lie groups due to Harish-Chandra \cite{HarishChandraplancherelcomplex}, see also Section 6.1 in \cite{Varadarajanbook}. 
Firstly, the characters of principal series representations are shown to be related to orbital integrals using Fourier transform. 
In a second step, orbital integrals on the group are transported to the Lie algebra. The final ingredient in the argument is 
the limit formula for orbital integrals on the Lie algebra, which in combination with the Weyl integration formula completes the proof. 

Trying to adapt this strategy to the quantum case seems difficult for various reasons. In fact, it is not even clear how to define a suitable notion of 
orbital integrals in this setting, and there is no good analogue of the Lie algebra. We proceed by explicitly writing down candidates for the 
Plancherel measure and Duflo-Moore operators instead, generalizing the ones in \cite{BRLorentz}. In order to verify that our choices are correct, we determine 
the characters of principal series representations and define a certain linear functional on the algebra of functions on the quantum group, which we 
call Plancherel functional. According to the Plancherel inversion formula it then suffices to show that the Plancherel 
functional agrees with the counit. 

For this, in turn, we use the BGG complex for quantized universal enveloping algebras, first described in \cite{ROSSO_BGGsln}, \cite{FEIGIN_FRENKEL_Toda}, \cite{MALIKOV_singularvectors}, and carefully studied by Heckenberger and Kolb in \cite{HeckenbergerKolbBGG}.  In its original form, the BGG complex is a resolution of a finite dimensional representation by direct sums of Verma modules, but here we will be interested in the geometric version, which is classically presented as a differential complex of sections of induced line bundles over a flag variety.  Our BGG complex is the resolution of the trivial representation by parabolically induced representations, and is obtained from the algebraic BGG complex via the category equivalence between category $ \O $ and the category of Harish-Chandra modules due to Joseph and Letzter \cite{JLVermaannihilator}, see also \cite{Josephbook}, \cite{VYcqg}. The key fact that allows us to compute the Plancherel functional is that its values can be identified 
with Lefschetz numbers of certain endomorphisms of the BGG complex. Since the BGG complex has almost trivial homology, an application of the 
Hopf trace formula completes the proof. 

Our result shows in particular that the Plancherel measure of complex semisimple quantum groups is supported on the space of unitary 
principal series representations, in analogy with the classical situation. This allows us to identify the reduced group $ C^* $-algebras 
of these quantum groups explicitly with certain continuous bundles of algebras of compact operators. As a consequence, one obtains a very transparent 
illustration of the deformation aspect in the operator algebraic approach to complex semisimple quantum groups, a feature which is not at all visible from 
the Drinfeld double construction.  

Let us now explain how the paper is organized. In Section \ref{secprelim} we collect some preliminaries on quantum groups and 
fix our notation. Section \ref{secrep} covers more specific background on complex semisimple quantum groups and their representations. 
We introduce our candidate Duflo-Moore operators for these quantum groups and compute the corresponding twisted characters of unitary principal 
series representations. In Section \ref{secplancherel} we recall the abstract Plancherel Theorem for locally compact quantum groups due to Desmedt. 
Section \ref{secplancherelcqg} contains our main result, that is, the Plancherel formula for complex semisimple quantum groups. As already indicated above, 
the proof involves the BGG complex, and we review the necessary background material along the way. 
In Section \ref{secproofslq2} we make some further comments and discuss a slightly different, more direct proof of the Plancherel formula in the simplest 
special case of the quantum Lorentz group. This argument is considerably shorter than the original proof by Buffenoir and Roche. 
Finally, in Section \ref{secreducedcstar} we apply the Plancherel formula to obtain an explicit description of the reduced group $ C^* $-algebras 
of arbitrary complex semisimple quantum groups. 

Let us conclude with some remarks on notation. The algebra of adjointable operators on a Hilbert space or Hilbert module $ \E $ is denoted 
by $ \LH(\E) $, and we write $ \KH(\E) $ for the algebra of compact operators. 
Depending on the context, the symbol $ \otimes $ denotes the algebraic tensor product over the 
complex numbers, the tensor product of Hilbert spaces, or the minimal tensor product of $ C^* $-algebras.

\section{Preliminaries} \label{secprelim}

In this section we review some background material on quantum groups and fix our notation. For more details 
we refer to \cite{CPbook}, \cite{KS}, \cite{KVLCQG}, \cite{VYcqg}. 

Throughout we assume that our definition parameter $ q $ is a strictly positive real number and $ q \neq 1 $. We write 
$$
[z]_q = \frac{q^z - q^{-z}}{q - q^{-1}}
$$
for the $ q $-number associated with $ z \in \mathbb{C} $ and use standard definitions and notation from $ q $-calculus. 

Let $ G $ be a simply connected complex semisimple Lie group with Lie algebra $ \mathfrak{g} $, 
and let $ \mathfrak{k} \subset \mathfrak{g} $ be the Lie algebra of a maximal compact subgroup $ K \subset G $. 
We fix a Cartan subalgebra $ \mathfrak{h} \subset \mathfrak{g} $ and a maximal torus $ T $ of $ K $ with Lie algebra $ \mathfrak{t} $ 
such that $ \mathfrak{t} \subset \mathfrak{h} $.
Let us denote by $ \simpleroots = \{\alpha_1, \dots, \alpha_N\} $ a set of simple roots for $ \mathfrak{g} $, 
and let $ (\;,\;) $ be the bilinear form on $ \mathfrak{h}^* $ obtained 
by rescaling the Killing form such that the shortest root $ \alpha $ of $ \mathfrak{g} $ satisfies $ (\alpha, \alpha) = 2 $. 
The simple coroots are given by $ \alpha_i^\vee = d_i^{-1} \alpha_i $ where $ d_i = (\alpha_i, \alpha_i)/2 $, 
and the entries of the Cartan matrix of $ \mathfrak{g} $ are $ a_{ij} = (\alpha_i^\vee, \alpha_j) $. 
We write $ \varpi_1, \dots, \varpi_N $ for the fundamental weights, defined by stipulating $ (\varpi_i, \alpha_j^\vee) = \delta_{ij} $.
Moreover we denote by $ \roots \subset \weights \subset \mathfrak{h}^* $ the root and weight lattices of $ \mathfrak{g} $, 
respectively. 
The set $ \weights^+ \subset \weights $ of dominant integral weights consists of all non-negative integer combinations of the fundamental 
weights. 

\begin{definition} \label{defuqg}
The quantized universal enveloping algebra $ U_q(\mathfrak{g}) $ is the complex algebra 
with generators $ K_\lambda $ for $ \lambda \in \weights $, and $ E_i, F_i $ for $ i = 1, \dots, N $, and the defining relations
\begin{align*}
K_0 &= 1 \\
K_\lambda K_\mu &= K_{\lambda + \mu} \\
K_\lambda E_j K_\lambda^{-1} &= q^{(\lambda, \alpha_j)} E_j \\
K_\lambda F_j K_\lambda^{-1} &= q^{-(\lambda, \alpha_j)} F_j \\ 
[E_i, F_j] &= \delta_{ij} \frac{K_i - K_i^{-1}}{q_i - q_i^{-1}} 
\end{align*} 
for all $ \lambda, \mu \in \weights $ and all $ i, j $, together with the quantum Serre relations 
\begin{align*}
&\sum_{k = 0}^{1 - a_{ij}} (-1)^k \begin{bmatrix} 1 - a_{ij} \\ k \end{bmatrix}_{q_i}
E_i^k E_j E_i^{1 - a_{ij} - k} = 0 \\ 
&\sum_{k = 0}^{1 - a_{ij}} (-1)^k \begin{bmatrix} 1 - a_{ij} \\ k \end{bmatrix}_{q_i}
F_i^k F_j F_i^{1 - a_{ij} - k} = 0. 
\end{align*}
In the above formulas we abbreviate $ K_i = K_{\alpha_i} $ for all simple roots, and we use the notation $ q_i = q^{d_i} $. 
\end{definition} 

We consider the Hopf algebra structure on $ U_q(\mathfrak{g}) $ determined by the 
comultiplication $ \hat{\Delta}: U_q(\mathfrak{g}) \rightarrow U_q(\mathfrak{g}) \otimes U_q(\mathfrak{g}) $ given by 
\begin{align*}
\hat{\Delta}(K_\lambda) &= K_\lambda \otimes K_\lambda \\
\hat{\Delta}(E_i) &= E_i \otimes K_i + 1 \otimes E_i \\
\hat{\Delta}(F_i) &= F_i \otimes 1 + K_i^{-1} \otimes F_i,  
\end{align*}
counit $ \hat{\varepsilon}: U_q(\mathfrak{g}) \rightarrow \mathbb{C} $ given by 
\begin{align*}
\hat{\varepsilon}(K_\lambda) = 1, \qquad \hat{\varepsilon}(E_j) = 0, \qquad \hat{\varepsilon}(F_j) = 0, 
\end{align*} 
and antipode $ \hat{S}: U_q(\mathfrak{g}) \rightarrow U_q(\mathfrak{g}) $ given by 
\begin{align*}
\hat{S}(K_\lambda) = K_{-\lambda}, \qquad \hat{S}(E_j) = -E_j K_j^{-1}, \qquad \hat{S}(F_j) = -K_j F_j
\end{align*}
on generators. We will use the Sweedler notation $ \hat{\Delta}(X) = X_{(1)} \otimes X_{(2)} $ 
for the comultiplication of $ U_q(\mathfrak{g}) $. 

We denote by $ U_q(\mathfrak{h}) $ the subalgebra of $ U_q(\mathfrak{g}) $ generated by the elements $ K_\lambda $ for $ \lambda \in \weights $, 
and let $ \mathfrak{h}^*_q $ be the space of all algebra characters $ U_q(\mathfrak{h}) \rightarrow \mathbb{C} $. 
Every such character is of the form $ \chi_\lambda(K_\mu) = q^{(\lambda, \mu)} $ for some $ \lambda \in \mathfrak{h}^* $,  
and if we write $ q = e^h $ and $ \hbar = \tfrac{h}{2\pi} $ we obtain an identification
$$ 
\mathfrak{h}^*_q = \mathfrak{h}^*/i\hbar^{-1}\roots^\vee
$$
in this way, where $ \roots^\vee $ is the coroot lattice.  

Let $ V $ be a left module over $ U_q(\mathfrak{g}) $. For $ \lambda \in \mathfrak{h}_q^* $ we define the weight space 
$$
V_\lambda = \{v \in V \mid K_\mu \cdot v = q^{(\mu, \lambda)} v \; \text{for all} \; \mu \in \weights \}. 
$$
We say that $ \lambda $ is a weight of $ V $ if $ V_\lambda $ is nonzero. 
A vector $ v \in V $ is said to have weight $ \lambda $ iff $ v \in V_\lambda $. A highest weight vector is a weight vector $ v $ such that 
$ E_i \cdot v = 0 $ for $ 1, \dots, N $. 
A module $ V $ over $ U_q(\mathfrak{g}) $ is called a weight module if it is the direct sum of its weight spaces $ V_\lambda $ 
for $ \lambda \in \mathfrak{h}_q^* $. 

The Verma module $ M(\lambda) $ is the universal weight module over $ U_q(\mathfrak{g}) $ generated by a highest weight vector $ v_\lambda $ 
of weight $ \lambda \in \mathfrak{h}^*_q $. As in the classical case it admits a unique irreducible quotient $ V(\lambda) $.  

We say that a weight module $ V $ is integrable if the operators $ E_i, F_j $ are locally nilpotent on $ V $ 
for all $ 1 \leq i,j \leq N $ and the weights of $ V $ are all contained in $ \weights \subset \mathfrak{h}^*_q $. 
Every finite dimensional weight module is completely reducible, and the irreducible integrable finite dimensional weight modules of $ U_q(\mathfrak{g}) $ are 
parametrized by their highest weights in $ \weights^+ $ as in the classical theory. If $ \mu \in \weights^+ $ we will 
write $ \pi_\mu: U_q(\mathfrak{g}) \rightarrow \End(V(\mu)) $ for the corresponding representation. The direct sum of the maps $ \pi_\mu $ induces 
an embedding $ \pi: U_q(\mathfrak{g}) \rightarrow \prod_{\mu \in \weights^+} \End(V(\mu)) $. 

The space of all matrix coefficients of finite dimensional integrable weight modules over $ U_q(\mathfrak{g}) $ is denoted by $ \Poly(G_q) $. 
It becomes a Hopf algebra with multiplication, comultiplication, counit and antipode in such a way that the 
canonical evaluation $ U_q(\mathfrak{g}) \times \Poly(G_q) \rightarrow \mathbb{C} $ is a skew-pairing, that is, we have 
\begin{align*}
(XY, f) &= (X, f_{(1)}) (Y, f_{(2)}), \qquad (X, fg) = (X_{(2)}, f) (X_{(1)}, g) 
\end{align*}
and 
\begin{align*}
(\hat{S}(X), f) &= (X, S^{-1}(f)), \qquad (\hat{S}^{-1}(X), f) = (X, S(f)) 
\end{align*} 
for $ X, Y \in U_q(\mathfrak{g}) $ and $ f, g \in \Poly(G_q) $. Here we use the Sweedler notation $ \Delta(f) = f_{(1)} \otimes f_{(2)} $ for 
the coproduct of $ f \in \Poly(G_q) $, and write $ S, \varepsilon $ for the antipode and counit of $ \Poly(G_q) $. 

Let us next discuss $ * $-structures. The quantized universal enveloping algebra $ U_q(\mathfrak{g}) $ is a Hopf $ * $-algebra with $ * $-structure given by 
\begin{align*}
E_i^* = K_i F_i, \qquad F_i^* = E_i K_i^{-1}, \qquad K_\lambda^* = K_\lambda.
\end{align*}
With the above $ * $-structure, $ U_q(\mathfrak{g}) $ should be viewed as the quantized universal enveloping algebra of the complexification 
of $ \mathfrak{k} $, and as such we shall write $ U_q^\mathbb{R}(\mathfrak{k}) $ for $ U_q(\mathfrak{g}) $ when we consider it as a Hopf $ * $-algebra. 
The representations $ V(\mu) $ for $ \mu \in \weights $ are $ * $-representations with respect to a uniquely determined inner product on $ V(\mu) $ 
for which the highest weight vector $ v_\mu $ has norm $ 1 $. 

Dually, we obtain a Hopf $ * $-algebra structure on $ \Poly(G_q) $ by stipulating
\begin{align*}
(X, f^*) &= \overline{(\hat{S}^{-1}(X)^*, f)} 
\end{align*} 
for $ f \in \Poly(G_q) $ and $ X \in U_q^\mathbb{R}(\mathfrak{k}) $. We will write $ \CF^\infty(K_q) $ for
$ \Poly(G_q) $ when we consider it as a Hopf $ * $-algebra in this way. The canonical bilinear pairing between $ U_q(\mathfrak{g}) $ and $ \Poly(G_q) $ then 
defines a skew-pairing of the Hopf $ * $-algebras $ U_q^\mathbb{R}(\mathfrak{k}) $ and $ \CF^\infty(K_q) $. 
The algebra $ \CF^\infty(K_q) $ can be viewed as a deformation of the Hopf $ * $-algebra of representative functions on the compact group $ K $. 
We will also write $ U_q^\mathbb{R}(\mathfrak{t}) $ for the Hopf $ * $-subalgebra of $ U_q^\mathbb{R}(\mathfrak{k}) $ 
with underlying algebra $ U_q(\mathfrak{h}) $. 

For each $ \mu \in \weights^+ $ we fix an orthonormal basis $ e_1^\mu, \dots, e_n^\mu $ of $ V(\mu) $ consisting of weight 
vectors, where $ n = \dim(V(\mu)) $. Then the formulas 
$$ 
(X, u^\mu_{ij}) = \bra e_i^\mu, \pi_\mu(X)(e_j^\mu) \ket = \bra e_i^\mu, X \cdot e_j^\mu \ket 
$$ 
define the corresponding matrix coefficients $ u^\mu_{ij} \in \CF^\infty(K_q) $, and we note that $ (u_{ij}^\mu)^* = S(u^\mu_{ji}) $. 
If $ \rho \in \weights $ denotes the half-sum of all positive roots then the quantum dimension of $ V(\mu) $ is defined by 
$$ 
\dim_q(V(\mu)) = \tr_{V(\mu)}(\pi_\mu(K_{2 \rho})) = \tr_{V(\mu)}(\pi_\mu(K_{-2 \rho})) = \sum_{j = 1}^n (K_{-2\rho}, u^\mu_{jj}), 
$$ 
where $ \tr_{V(\mu)} \in \CF^\infty(K_q) $ denotes the natural trace on $ V(\mu) $. 
If we write $ \phi $ for the Haar state of $ \CF^\infty(K_q) $, then the Schur orthogonality relations are
$$
\phi(u^\beta_{ij} S(u^\gamma_{lk})) = \delta_{\beta \gamma} \delta_{ik} \, \frac{(K_{-2 \rho}, u^\beta_{jl})}{\dim_q(V(\beta))}, 
\qquad 
\phi(S(u^\beta_{ji}) u^\gamma_{kl}) = \delta_{\beta \gamma} \delta_{jl} \, \frac{(K_{2 \rho}, u^\beta_{ik})}{\dim_q(V_\beta)}, 
$$
compare for instance Chapter 11 in \cite{KS}. These relations imply the modular property 
$$ 
\phi(fg) = (K_{2 \rho}, g_{(1)} g_{(3)}) \phi(g_{(2)} f) 
$$ 
for all $ f, g \in \CF^\infty(K_q) $. 

The Hopf $ * $-algebra $ \CF^\infty(K_q) $ is an algebraic quantum group in the sense of Van Daele \cite{vDadvances}, that is, a multiplier Hopf $ * $-algebra
with nonzero positive left invariant Haar functional. We write $ \DF(K_q) $ for the dual algebraic quantum group. 
Explicitly, the dual is given by the algebraic direct sum 
$$
\DF(K_q) = \text{alg-}\bigoplus_{\mu \in \weights^+} \KH(V(\mu))  
$$
with the $ * $-structure arising from the $ C^* $-algebras $ \KH(V(\mu)) = \End(V(\mu)) $. We denote by $ p_\eta $ the central projection in $ \DF(K_q) $ 
corresponding to the matrix block $ \KH(V(\eta)) $ for $ \eta \in \weights^+ $. 

There exists a unique bilinear pairing $ \DF(K_q) \times \CF^\infty(K_q) \rightarrow \mathbb{C} $ such that 
\begin{align*}
(xy, f) &= (x, f_{(1)}) (y, f_{(2)}), \qquad (x, fg) = (x_{(2)}, f) (x_{(1)}, g) 
\end{align*}
and 
\begin{align*}
(\hat{S}(x), f) &= (x, S^{-1}(f)), \qquad (\hat{S}^{-1}(x), f) = (x, S(f)) 
\end{align*} 
for $ f, g \in \CF^\infty(K_q) $ and $ x, y \in \DF(K_q) $. 
The compatibility with the $ * $-structures is given by 
\begin{align*}
(x, f^*) &= \overline{(\hat{S}^{-1}(x)^*, f)}, \qquad (x^*, f) = \overline{(x, S(f)^*)}. 
\end{align*}
Positive left and right Haar functionals for $ \DF(K_q) $ are given by 
\begin{align*}
\hat{\phi}(x) &= \sum_{\mu \in \weights^+} \dim_q(V(\mu)) \tr(K_{2\rho} p_\mu x), \qquad 
\hat{\psi}(x) = \sum_{\mu \in \weights^+} \dim_q(V(\mu)) \tr(K_{-2\rho} p_\mu x), 
\end{align*}
respectively. 

Let us write 
$$
\M(\DF(K_q)) = \text{alg-}\prod_{\mu \in \weights^+} \KH(V(\mu))  
$$  
for the algebraic multiplier algebra of $ \DF(K_q) $. 
The pairing between $ \DF(K_q) $ and $ \CF^\infty(K_q) $ extends uniquely to a bilinear pairing between $ \M(\DF(K_q)) $ and $ \CF^\infty(K_q) $. 
If we consider the canonical embedding $ U_q^\mathbb{R}(\mathfrak{k}) \subset \M(\DF(K_q)) $, then this is compatible with 
our original pairing between $ U_q^\mathbb{R}(\mathfrak{k}) $ and $ \CF^\infty(K_q) $. 

By Pontrjagin duality, we can also view $ \DF(K_q) $ as function algebra of the dual algebraic quantum group $ \hat{K}_q $, and $ \CF^\infty(K_q) $ 
as its dual. However, when one flips the roles of the two algebras one has to be slightly careful.  
In particular, the natural pairing $ \CF^\infty(K_q) \times \DF(K_q) \rightarrow \mathbb{C} $ is defined by 
$$
(f, x) = (\hat{S}(x), f) = (x, S^{-1}(f)) 
$$
for $ f \in \CF^\infty(K_q) $ and $ x \in \DF(K_q) $. The antipode is needed in order to obtain the skew-pairing property and the correct 
behaviour with respect to the $ * $-structures on both sides. We emphasize that, with these conventions, we have $ (f,x) \neq (x, f) $ in general. 

Given the basis of matrix coefficients $ u^\mu_{ij} $ in $ \CF^\infty(K_q) $ as above we obtain a dual linear basis of matrix 
units $ \omega^\mu_{ij} $ of $ \DF(K_q) $ satisfying 
$$
(\omega^\mu_{ij}, u^\nu_{kl}) = \delta_{\mu \nu} \delta_{ik} \delta_{jl}. 
$$
The fundamental multiplicative unitary of the quantum group $ K_q $ is the algebraic multiplier of $ \CF^\infty(K_q) \otimes \DF(K_q) $ given by 
$$
W = \sum_{\mu \in \weights^+} \sum_{i,j = 1}^{\dim(V(\mu))} u^\mu_{ij} \otimes \omega^\mu_{ij}, 
$$
and we have the formula 
$$
W^{-1} = (S \otimes \id)(W) = (\id \otimes \hat{S}^{-1})(W)
$$
for its inverse. 

With these preparations in place, let us now discuss the main object of study in this paper, namely the Drinfeld double $ G_q = K_q \bowtie \hat{K}_q $. 
By definition, this is the algebraic quantum group given by the $ * $-algebra
$$
\CF^\infty_c(G_q) = \CF^\infty(K_q) \otimes \DF(K_q),  
$$
with comultiplication 
$$
\Delta_{G_q} = (\id \otimes \sigma \otimes \id)(\id \otimes \ad(W) \otimes \id)(\Delta \otimes \hat{\Delta}), 
$$
counit 
$$ 
\varepsilon_{G_q} = \varepsilon \otimes \hat{\varepsilon}, 
$$
and antipode 
$$
S_{G_q}(f \otimes x) = W^{-1} (S(f) \otimes \hat{S}(x)) W = (S \otimes \hat{S})(W (f \otimes x) W^{-1}). 
$$
Here $ W \in \M(\CF^\infty(K_q) \otimes \DF(K_q)) $ denotes the multiplicative unitary from above. 
A positive left and right invariant Haar functional for $ \CF^\infty_c(G_q) $ is given by 
$$
\phi_{G_q}(f \otimes x) = \phi(f) \hat{\psi}(x), 
$$
compare \cite{PWlorentz}. 

Dually, we obtain the convolution algebra $ \DF(G_q) = \DF(K_q) \bowtie \CF^\infty(K_q) $, which has $ \DF(K_q) \otimes \CF^\infty(K_q) $ 
as its underlying vector space, equipped with the tensor product comultiplication and the multiplication 
\begin{align*}
(x \bowtie f)(y \bowtie g) &= x (y_{(1)}, f_{(1)}) y_{(2)} \bowtie f_{(2)} (\hat{S}(y_{(3)}), f_{(3)}) g.  
\end{align*}
The $ * $-structure of $ \DF(G_q) $ is defined in such a way that the natural inclusion homomorphisms $ \DF(K_q) \rightarrow \DF(G_q) $ and 
$ \CF^\infty(K_q) \rightarrow \M(\DF(G_q)) $ are $ * $-homomorphisms. 
By construction these maps are compatible with the comultiplications,
and we record the formula
$$
\hat{\varepsilon}_{G_q}(x \bowtie f) = \hat{\varepsilon}(x) \varepsilon(f)
$$
for the counit of $ \DF(G_q) $. 
We shall work with the skew-pairing 
$$
(y \bowtie g, f \otimes x) = (y, f)(g, x)
$$
between $ \DF(G_q) $ and $ \CF^\infty_c(G_q) $, and we remark that this is compatible with the $ * $-structures defined above.  
For more details on the properties of this pairing we refer to Chapter 4 of \cite{VYcqg}. 

Both $ * $-algebras $ \CF^\infty_c(G_q) $ and $ \DF(G_q) $ admit universal $ C^* $-completions, 
which will be denoted by $ C_0(G_q) $ and $ C^*_\max(G_q) $, respectively. 
By definition, a unitary representation of $ G_q $ on a Hilbert space $ \H $ 
is a nondegenerate $ * $-homomorphism $ \pi: C^*_\max(G_q) \rightarrow \LH(\H) $. A basic example is the left regular representation 
of $ G_q $, which is obtained from the canonical $ * $-homomorphism $ C^*_\max(G) \rightarrow \LH(L^2(G_q)) $. 
Here $ L^2(G_q) $ is the GNS-construction of the left Haar functional of $ G_q $. 
The image of $ C^*_\max(G) $ under the regular representation is the reduced group $ C^* $-algebra $ C^*_\red(G_q) \subset \LH(L^2(G_q)) $. 

Let us finally recall that the group algebra $ \DF(G_q) $ can be identified with the vector space $ \CF^\infty_c(G_q) $ equipped with the convolution product 
$$
f * g = \phi_{G_q}(S_{G_q}^{-1}(g_{(1)}) f) g_{(2)} 
= \phi_{G_q}(S_{G_q}^{-1}(g) f_{(2)}) f_{(1)} 
$$
by Fourier transform, see Section 4.2.2 in \cite{VYcqg}. More precisely, the linear map $ \F: \CF^\infty_c(G_q) \rightarrow \DF(G_q) $ determined by 
$$
(\F(f), h) = \phi_{G_q}(hf)
$$
is a linear isomorphism which identifies $ \CF^\infty_c(G_q) $, equipped with the convolution product, with $ \DF(G_q) $ as algebras. 
It becomes an isomorphism of $ * $-algebras if we consider the $ * $-structure defined by $ f^* = \F^{-1}(\F(f)^*) $ for $ f \in \CF^\infty_c(G_q) $, 
not to be confused with the $ * $-structure underlying $ \CF^\infty_c(G_q) $. We will mainly work with this description of the group algebra $ \DF(G_q) $ 
in our calculations below.

\section{Representation theory of complex quantum groups} \label{secrep}

In this section we review some central facts regarding the representation theory of complex quantum groups. For the proofs of these 
results as well as further background we refer to \cite{VYcqg}. Throughout we assume that $ q = e^h $ is a strictly positive real number and $ q \neq 1 $. 

Let $ \mu \in \weights $. Then we define the space of sections $ \Gamma(\E_\mu) \subset \CF^\infty(K_q) $ of the induced vector bundle 
$ \E_\mu $ corresponding to $ \mu $ to be the subspace of $ \CF^\infty(K_q) $ of weight $ \mu $ with respect to the $ U_q^\mathbb{R}(\mathfrak{k}) $-module 
structure 
$$
X \hit \xi =\xi_{(1)} (X, \xi_{(2)}).  
$$
Equivalently, we have 
$$
\Gamma(\E_\mu) = \{\xi \in \CF^\infty(K_q) \mid (\id \otimes \pi_T) \Delta(\xi) = \xi \otimes e^\mu \}, 
$$ 
where $ \pi_T: \CF^\infty(K_q) \rightarrow \CF^\infty(T) $ is the canonical projection homomorphism and $ e^\mu \in \CF^\infty(T) $ 
is the generator corresponding to the weight $ \mu $. We note that $ \CF^\infty(T) $ is the quotient of $ \CF^\infty(K_q) $ corresponding to the 
Hopf $ * $-subalgebra $ U_q^\mathbb{R}(\mathfrak{t}) \subset U_q^\mathbb{R}(\mathfrak{k}) $. 

For $ \lambda \in \mathfrak{h}^*_q $ we define the twisted left adjoint action of $ \CF^\infty(K_q) $ on $ \Gamma(\E_\mu) $ by 
\begin{equation*} 
f \cdot \xi = f_{(1)} \xi S(f_{(3)}) \,(K_{2 \rho + \lambda}, f_{(2)}).
\end{equation*}
Together with the left coaction $ \Gamma(\E_\mu) \rightarrow \CF^\infty(K_q) \otimes \Gamma(\E_\mu) $ given by comultiplication  
this turns $ \Gamma(\E_\mu) $ into a Yetter-Drinfeld module. 
We will frequently switch from the left coaction on $ \Gamma(\E_\mu) $ to the left $ \DF(K_q) $-module structure given by 
$$
x \cdot \xi = (\hat{S}(x), \xi_{(1)}) \xi_{(2)} 
$$
for $ x \in \DF(K_q) $. Combining this with the action of $ \CF^\infty(K_q) $ from above makes $ \Gamma(\E_\mu) $ into a $ \DF(G_q) $-module, 
which we denote by $ \Gamma(\E_{\mu, \lambda}) $ and refer to as the \emph{principal series module} 
with parameter $ (\mu, \lambda) \in \weights \times \mathfrak{h}_q^* $.  

Let us write $ \mathfrak{t}^*_q = \mathfrak{t}^*/i\hbar^{-1}\roots^\vee \subset \mathfrak{h}^*_q $, 
where, by slight abuse of notation, we view the dual space $ \mathfrak{t}^* = \Hom_\mathbb{R}(\mathfrak{t}, \mathbb{R}) $ as real vector subspace 
of $ \mathfrak{h}^* = \Hom_{\mathbb{C}}(\mathfrak{h}, \mathbb{C}) \cong \Hom_\mathbb{R}(\mathfrak{t}, \mathbb{C}) $. 
It will also be convenient to consider $ \mathfrak{t}^* = i \mathfrak{a}^* $ and work with $ \mathfrak{a}^*_q = \mathfrak{a}^*/\hbar^{-1}\roots^\vee $, 
so that $ i \mathfrak{a}^*_q \subset \mathfrak{h}^*_q $ can be identified with $ \mathfrak{t}^*_q $. 

For $ \lambda \in \mathfrak{t}^*_q $, or equivalently $ \lambda = i \nu $ for $ \nu \in \mathfrak{a}_q^* $, the Yetter-Drinfeld 
module $ \Gamma(\E_{\mu, \lambda}) \subset \CF^\infty(K_q) $ is unitary for the standard scalar 
product on $ \CF^\infty(K_q) $. In particular, we obtain a corresponding 
nondegenerate $ * $-representation $ \pi_{\mu, i \nu}: C^*_\max(G_q) \rightarrow \LH(\H_{\mu, i\nu}) $, 
where $ \H_{\mu, i\nu} \subset L^2(K_q) $ is the Hilbert space completion of $ \Gamma(\E_{\mu, i\nu}) $. 

\begin{definition} 
The unitary representations of $ G_q $ on $ \H_{\mu, i \nu} $ for $ (\mu, \nu) \in \weights \times \mathfrak{a}^*_q $ as above 
are called unitary principal series representations. 
\end{definition}

For proofs of the following results we refer to Chapter 6 of \cite{VYcqg}. 

\begin{theorem} \label{thmprincipalirreducible}
For all $ (\mu, \nu) \in \weights \times \mathfrak{a}^*_q $ the unitary principal series representation $ \H_{\mu, i\nu} $ is an 
irreducible representation of $ G_q $. 
\end{theorem} 

The Weyl group $ W $ acts on the parameter space $ \weights \times \mathfrak{a}^*_q $ by 
$$
w (\mu, \nu) = (w \mu, w \nu). 
$$
The following result describes the isomorphisms between unitary principal series representations in the quantum case. 

\begin{theorem} \label{thmweylgroupactionirreds}
Let $ (\mu, \nu), (\mu', \nu') \in \weights \times \mathfrak{a}^*_q $. Then $ \H_{\mu, i \nu} $ and $ \H_{\mu', i \nu'} $ are equivalent 
representations of $ G_q $ iff $ (\mu, \nu), (\mu', \nu') $ are in the same Weyl group orbit, that is, iff 
$$ 
(\mu', \nu') = (w \mu, w \nu) 
$$ 
for some $ w \in W $. 
\end{theorem} 

In the remainder of this section we shall study the characters of unitary principal series representations. 
Fix $ (\mu, \nu) \in \weights \times \mathfrak{a}_q^* $, and let $ f $ be an element of the convolution algebra $ \CF^\infty_c(G) $. 
We will write $ \pi_{\mu, i\nu}(f) $ for the corresponding operator on $ \H_{\mu, i\nu} $, by identifying $ f $ with 
an element of $ \DF(G_q) \subset C^*_\max(G_q) $ as explained at the end of Section \ref{secprelim}. 
It follows from the structure of induced bundles that $ \pi_{\mu, i\nu}(f) $ is a 
finite rank operator. In particular, the operator $ \pi_{\mu, i\nu}(f) $ is trace-class, and we shall be interested in 
finding an explicit formula for a certain twisted version of its operator trace. 

Let us recall that the homogeneous vector bundle $ \E_\mu $ can be described both using sections of an associated vector bundle over $ K_q/T_q $ 
and as sections of an associated vector bundle over $ G_q/B_q $, where $ B_q $ denotes the quantum analogue of the minimal parabolic 
subgroup $ B \subset G $. The latter is defined as the relative Drinfeld double $ T \bowtie \hat{K}_q $, so 
that $ \CF^\infty_c(B_q) = \CF^\infty(T) \otimes \DF(K_q) $, with a suitable twisted comultiplication. 

Our definition of $ \Gamma(\E_{\mu, i\nu}) = \Gamma(\E_\mu) $ above was phrased in the compact picture, namely
$$ 
\Gamma(\E_{\mu, i\nu}) = \{\xi \in \CF^\infty(K_q) \mid (\id \otimes \pi_T) \Delta(\xi) = \xi \otimes e^\mu \}. 
$$ 
In the noncompact picture, we consider instead elements $ \sigma $ of the algebraic multiplier algebra $ \CF^\infty(G_q) $ of $ \CF^\infty_c(G_q) $ such that 
$$ 
(\id \otimes \pi_{B_q}) \Delta_{G_q}(\sigma) = \sigma \otimes (e^\mu \otimes K_{2\rho + i\nu}).
$$ 
Here $ K_{2\rho + i\nu} $ is viewed as multiplier of $ \DF(K_q) $ inside $ \CF^\infty(G_q) = \M(\CF^\infty_c(G_q)) $, 
and $ \pi_{B_q}: \CF^\infty_c(G_q) \rightarrow \CF^\infty_c(B_q) $ is the canonical projection. 
If $ \xi \in \Gamma(\E_{\mu, i\nu}) $ then the corresponding element $ \ext(\xi) $ in the noncompact picture is given by 
$$ 
\ext(\xi) = \xi \otimes K_{2\rho + i\nu} \in \CF^\infty(G_q). 
$$
Conversely, if $ \sigma \in \CF^\infty(G_q) $ satisfies the invariance condition in the noncompact picture then 
$$ 
\res(\sigma) = (\id \otimes \hat{\varepsilon})(\sigma) 
$$ 
is contained in $ \Gamma(\E_{\mu, i\nu}) $, and the maps $ \ext $ and $ \res $ are inverse to each other \cite[Lemma 6.18]{VYcqg}. 

Recall that we may identify the group algebra $ \DF(G_q) $ with $ \CF^\infty_c(G_q) $ using Fourier transform, where the latter is equipped with convolution.   
With this in mind, the action of $ f \in \CF^\infty_c(G_q) $ on $ \xi \in \Gamma(\E_{\mu, i\nu}) $ is given by 
\begin{align*}
\pi_{\mu, i \nu}(f)(\xi) 
&= \phi_{G_q}(S_{G_q}^{-1}(\ext(\xi)) f_{(2)}) \res(f_{(1)}). 
\end{align*}
In particular, if $ f = a \otimes t $ for $ a \in \CF^\infty(K_q) $ and $ t \in \DF(K_q) $, then using 
$$
S_{G_q}^{-1}(\ext(\xi)) = S^{-1}_{G_q}(\xi \otimes K_{2\rho + i\nu}) = W^{-1}(S^{-1}(\xi) \otimes K_{-2\rho - i\nu}) W 
$$
we obtain 
\begin{align*}
\pi_{\mu, i \nu}(f)(\xi) &=
\sum_{\nu,\eta \in \weights^+} \sum_{m,n,r,s} 
\phi(S(u^\nu_{nm}) S^{-1}(\xi) u^\eta_{rs} a_{(2)}) \hat{\psi}(\omega^\nu_{nm} K_{-2\rho - i\nu} \omega^\eta_{rs} t) a_{(1)} \\
&= \sum_{\eta \in \weights^+} \sum_{n,r,s} 
\phi(S(u^\eta_{nr}) S^{-1}(\xi) u^\eta_{rs} a_{(2)}) \hat{\psi}(\omega^\eta_{nr} K_{-2\rho - i\nu} \omega^\eta_{rs} t) a_{(1)},  
\end{align*}
taking into account that the operator $ K_{-2\rho - i\nu} $ is diagonal in our chosen basis of the representation $ V(\eta) $. 
For $ f = u^\beta_{ij} \otimes \omega^\gamma_{kl} $ this formula reduces to 
\begin{align*}
\pi_{\mu, i \nu}(f)(\xi) &=
\sum_{m, r} \dim_q(V(\gamma))
\phi(S(u^\gamma_{lr}) S^{-1}(\xi) u^\gamma_{rk} u^\beta_{mj}) q^{(-2\rho, \varepsilon_l + \varepsilon_r)} q^{-(i\nu, \varepsilon_r)} u^\beta_{im} \\
&= \sum_{\varepsilon_m = \mu} \sum_r \dim_q(V(\gamma))
\phi(S(u^\gamma_{lr}) S^{-1}(\xi) u^\gamma_{rk} u^\beta_{mj}) q^{(-2\rho, \varepsilon_l + \varepsilon_r)} q^{-(i\nu, \varepsilon_r)} u^\beta_{im}. 
\end{align*}
Here we write $ \varepsilon_m $ for the weight of the basis vector $ e^\beta_m $ in the definition of the matrix coefficient $ u^\beta_{im} $, 
and note that for $ \xi \in \Gamma(\E_{\mu, i\nu}) $ only terms with $ \varepsilon_m = \mu $ give nonzero contributions in the expression on the right hand side 
by weight considerations. 

Let us now introduce certain operators which will turn out to be the Duflo-Moore operators for $ G_q $, 
compare Theorem \ref{thmcqgplancherel} below. 

\begin{definition} \label{defduflomoore}
For $ (\mu, \nu) \in \weights \times \mathfrak{a}^*_q $ we let $ D_{\mu, i\nu} $ be the unbounded linear operator 
in $ \H_{\mu, i\nu} $ given by 
$$
D_{\mu, i \nu}(\xi) = \pi_{\mu, i \nu}(K_{-\rho} \bowtie 1)(\xi) = (K_{\rho}, \xi_{(1)}) \xi_{(2)} 
$$
for $ \xi \in \Gamma(\E_{\mu, i \nu}) $. 
\end{definition}

Here $ K_{-\rho} \bowtie 1 $ is viewed as a multiplier of $ \DF(G_q) $ in the obvious way, and we use that the representation $ \pi_{\mu, i \nu} $ 
extends naturally to a representation of $ \M(\DF(G_q)) $ on $ \Gamma(\E_{\mu, i \nu}) $. 

Using the Peter-Weyl decomposition of $ \H_{\mu, i \nu} \subset L^2(K_q) $, it is straightforward to check that the formula in 
Definition \ref{defduflomoore} uniquely determines an unbounded strictly positive self-adjoint operator in $ \H_{\mu, i \nu} $, which will 
again be denoted by $ D_{\mu, i\nu} $. We observe that 
$$
D_{\mu, i \nu}^{-2}(\xi) = \pi_{\mu, i \nu}(K_{2\rho} \bowtie 1)(\xi) = (K_{-2\rho}, \xi_{(1)}) \xi_{(2)} 
$$
for $ \xi \in \Gamma(\E_{\mu, i \nu}) $. 

In the sequel we shall again tacitly identify the group algebra $ \DF(G_q) $ with $ \CF^\infty_c(G_q) $ equipped with convolution. 
With this notational convention in mind, we remark that it is straightforward to check that $ \pi_{\mu, i \nu}(f) D_{\mu, i \nu}^{-2} $ defines a finite 
rank operator on $ \H_{\mu, i \nu} $ for all $ f \in \CF^\infty_c(G_q) $. The operator trace $ \tr(\pi_{\mu, i \nu}(f) D_{\mu, i \nu}^{-2}) $ will be referred 
to as the \emph{twisted character} of $ \pi_{\mu, i \nu}(f) $. 

\begin{prop} \label{characterformula}
Let $ (\mu, \nu) \in \weights \times \mathfrak{a}^*_q $ and $ f = u^\beta_{ij} \otimes \omega^\gamma_{kl} \in \CF^\infty_c(G) $. Then the 
twisted character of $ \pi_{\mu, i \nu}(f) $ is given by 
\begin{align*}
\tr(&\pi_{\mu, i \nu}(f) D_{\mu, i \nu}^{-2}) \\
&= q^{(-2 \rho, \mu)} \dim_q(V(\gamma)) \sum_{\varepsilon_m = \mu} \sum_{r} 
\phi(u^\gamma_{lr} S^{-1}(u^\beta_{mj}) S^{-1}(u^\gamma_{rk}) u^\beta_{im}) q^{-(i \nu, \varepsilon_r)}. 
\end{align*}
\end{prop} 

\proof Using the antipode relation $ \hat{S}(X) = K_{2\rho} \hat{S}^{-1}(X) K_{-2\rho} $ we obtain 
\begin{align*}
(X, S^{-1}(D_{\mu, i \nu}^{-2}(\xi))) &= (\hat{S}(X), (K_{-2\rho}, \xi_{(1)}) \xi_{(2)}) \\
&= (K_{-2\rho} \hat{S}(X), \xi) \\
&= (\hat{S}^{-1}(X) K_{-2 \rho}, \xi) \\
&= q^{(-2 \rho, \mu)} (X, S(\xi)) 
\end{align*} 
for $ \xi \in \Gamma(\E_{\mu, i \nu}) $ and $ X \in U_q^\mathbb{R}(\mathfrak{k}) $. 
Inserting this into the formula for $ \pi_{\mu, i \nu}(f)(\xi) $ obtained above and applying the operator trace yields 
\begin{align*}
\tr(&\pi_{\mu, i \nu}(f) D_{\mu, i \nu}^{-2}) \\
&= q^{(-2 \rho, \mu)} \dim_q(V(\gamma)) \sum_{\varepsilon_m = \mu} \sum_r q^{(-2 \rho, \varepsilon_l + \varepsilon_r)} 
\phi(S(u^\gamma_{lr}) S(u^\beta_{im}) u^\gamma_{rk} u^\beta_{mj})  q^{-(i \nu, \varepsilon_r)} \\
&= q^{(-2 \rho, \mu)} \dim_q(V(\gamma)) \sum_{\varepsilon_m = \mu} \sum_r q^{(-2 \rho, \varepsilon_l + \varepsilon_r)} 
\phi(S^{-1}(u^\beta_{mj}) S^{-1}(u^\gamma_{rk}) u^\beta_{im} u^\gamma_{lr}) q^{-(i \nu, \varepsilon_r)} \\
&= q^{(-2 \rho, \mu)} \dim_q(V(\gamma)) \sum_{\varepsilon_m = \mu} \sum_{r} 
\phi(u^\gamma_{lr} S^{-1}(u^\beta_{mj}) S^{-1}(u^\gamma_{rk}) u^\beta_{im}) q^{-(i \nu, \varepsilon_r)}, 
\end{align*}
using invariance under the antipode and the modular property of the Haar functional. \qed

\section{The abstract Plancherel Theorem} \label{secplancherel}

In this section we review the abstract Plancherel theorem for locally compact quantum groups due to 
Desmedt. We refer to \cite{Desmedtthesis} for further information. 

Let us say that a locally compact quantum group $ G $ is second countable if $ L^2(G) $ is a separable Hilbert space. By definition, $ G $ is type I 
if the group $ C^* $-algebra $ C^*_\max(G) $ is type I. We write $ \Irr(G) $ for the space of equivalence classes of irreducible unitary representations 
of $ G $ with the Fell topology. If $ G $ is type I then the space $ \Irr(G) $ is a standard Borel space with the
Borel structure coming from the Fell topology, see Section 4.6 in \cite{Dixmiercstar}. 

Given $ \lambda \in \Irr(G) $ let us write $ \H_\lambda $ for the underlying Hilbert space of a representative of $ \lambda $ 
and $ \pi_\lambda: C^*_\max(G) \rightarrow \LH(\H_\lambda) $ for the corresponding nondegenerate $ * $-homomorphism. 
We shall also write $ HS(\H_\lambda) = \H_\lambda \otimes \overline{\H}_\lambda $ for the space of Hilbert-Schmidt operators on $ \H_\lambda $, 
which we consider as a representation of $ G_q $ with the action on the first tensor factor. 

The following statement is then a condensed version of Theorem 3.4.1 in \cite{Desmedtthesis}. 

\begin{theorem}[Plancherel Theorem] \label{thmabstractplancherel} 
Let $ G $ be a second countable locally compact quantum group of type I. Then there exists a standard measure 
$ m $ on $ \Irr(G) $, a measurable field of Hilbert spaces $ (\H_\lambda)_{\lambda \in \Irr(G)} $, a measurable 
field $ (D_\lambda)_{\lambda \in \Irr(G)} $ of self-adjoint strictly positive operators for $ (\H_\lambda)_{\lambda \in \Irr(G)} $, 
and an isometric $ G $-equivariant isomorphism 
$$ 
Q: L^2(G) \rightarrow \int^\oplus_{\Irr(G)} HS(\H_\lambda) dm(\lambda), 
$$ 
given by 
$$ 
Q(\hat{\Lambda}(x)) = \int^\oplus_{\Irr(G)} \pi_\lambda(x) D_\lambda^{-1} dm(\lambda) 
$$ 
on a certain dense subspace of $ L^2(G) \cap L^1(G) $. The Plancherel measure is unique up to equivalence, 
more precisely, the family of Duflo-Moore operators $ (D_\lambda)_{\lambda \in \Irr(G)} $ combined with $ m $ are unique up to mutual rescaling. 
\end{theorem} 

Here we write $ L^1(G) $ for the predual of the von Neumann algebra $ L^\infty(G) $ associated to $ G $, and 
$ \hat{\Lambda}: \N_{\hat{\phi}} \rightarrow L^2(G) $ denotes the GNS-map for dual left Haar weight $ \hat{\phi} $ of $ G $. 

If $ G $ is an algebraic quantum group in the sense of van Daele \cite{vDadvances} then the initial domain of definition of the map $ Q $ 
in Theorem \ref{thmabstractplancherel} contains the space $ \CF^\infty_c(G) $. 
For computational purposes the following version of the Plancherel inversion formula will be useful for us.

\begin{lemma} \label{plancherelformula}
Let $ G $ be a second countable algebraic quantum group of type I. Let $ m $ be a standard measure on $ \Irr(G) $ and 
$ (D_\lambda)_{\lambda \in \Irr(G)} $ a measurable field of self-adjoint strictly positive operators on a measurable field of Hilbert 
spaces $ (\H_\lambda)_{\lambda \in \Irr(G)} $. \\
Then $ m $ is the Plancherel measure with Duflo-Moore 
operators $ (D_\lambda)_{\lambda \in \Irr(G)} $ iff for all $ f \in \CF^\infty_c(G) $ the operator $ \pi_\lambda(f) D_\lambda^{-2} $ is trace-class 
for almost all $ \lambda \in \Irr(G) $ and the Plancherel inversion formula
$$
\varepsilon(f) = \int_{\Irr(G)} \tr(\pi_\lambda(f) D_\lambda^{-2}) dm(\lambda)
$$
holds. Here $ \tr $ denotes the operator trace. 
\end{lemma}

\proof We note again that using Fourier transform we tacitly identify the group algebra $ \DF(G) $ of $ G $ with $ \CF^\infty_c(G) $, the 
latter being equipped with the convolution product $ g * h = \phi(S^{-1}(h_{(1)}) g) h_{(2)} $ and the $ * $-structure inherited from $ \DF(G_q) $. 

Assume first that $ m $ is the Plancherel measure with corresponding Duflo-Moore operators $ (D_\lambda)_{\lambda \in \Irr(G)} $. 
If $ g, h \in \CF^\infty_c(G) $ then by Theorem \ref{thmabstractplancherel} we have 
\begin{align*}
\bra \Lambda(g), \Lambda(h) \ket &= \int_{\Irr(G)} \tr(D_\lambda^{-1} \pi_\lambda(g)^* \pi_\lambda(h) D_\lambda^{-1}) dm(\lambda) \\
&= \int_{\Irr(G)} \tr(\pi_\lambda(g^* * h) D_\lambda^{-2}) dm(\lambda).
\end{align*}
Moreover 
$$ 
\varepsilon(g^* * h) = \hat{\phi}(\F(g^* * h)) = \hat{\phi}(\F(g)^* \F(h)) = \phi(g^* h) = \bra \Lambda(g), \Lambda(h) \ket 
$$ 
by properties of the Fourier transform $ \F $, using that $ \hat{\phi}(\F(f)) = \varepsilon(f) $ for any $ f \in \CF^\infty_c(G) $, and 
keeping in mind the different $ * $-structures associated with convolution and multiplication in $ \CF^\infty_c(G) $. 
In other words, both sides of the Plancherel inversion formula agree on $ f = g^* * h $.
Since elements of this form span $ \CF^\infty_c(G) $ linearly we see that the Plancherel inversion formula holds for all $ f \in \CF^\infty_c(G) $. 

Conversely, if the Plancherel inversion formula holds for all $ f \in \CF^\infty_c(G) $, then reversing the previous argument shows that the 
formula for the map $ Q $ in Theorem \ref{thmabstractplancherel} defines an isometric linear map on the dense subspace $ \CF^\infty_c(G) \subset L^2(G) $. 
From this it follows that the measure $ m $ and the operators $ (D_\lambda)_{\lambda \in \Irr(G)} $ satisfy the 
properties listed in Theorem \ref{thmabstractplancherel}, and therefore are the Plancherel measure with corresponding Duflo-Moore operators. \qed 

Theorem \ref{thmabstractplancherel} provides a complete description of the regular representation of the quantum group $ G $ at an abstract level. 
A key problem in harmonic analysis is to compute the Plancherel measure and corresponding Duflo-Moore operators concretely, given 
a parametrization of the space of irreducible representations.

\section{The Plancherel formula for complex quantum groups} \label{secplancherelcqg}

In this section we given an explicit description of the Plancherel formula for complex semisimple quantum groups. We remark that these quantum groups 
are indeed type I, see Chapter 6 in \cite{VYcqg}, so that Theorem \ref{thmabstractplancherel} applies. 
As in the classical case, our result shows in particular that the support of the Plancherel measure for a complex semisimple 
quantum group is the space of unitary principal series representations. 

\begin{theorem} \label{thmcqgplancherel}
Let $ G_q $ be a complex semisimple quantum group. Moreover let $ \H = (\H_{\mu, i\nu})_{\mu, \nu} $ 
be the Hilbert space bundle over $ \weights \times \mathfrak{a}_q^* $ of unitary principal series representations of $ G_q $. Then there is 
a $ G $-equivariant unitary isomorphism 
$$ 
Q: L^2(G_q) \rightarrow \bigoplus_{\mu \in \weights} \int^\oplus_{\nu \in \mathfrak{a}_q^*} HS(\H_{\mu, i\nu}) dm_\mu(\nu)
$$ 
given by 
$$ 
Q(\hat{\Lambda}(x)) = \sum_{\mu \in \weights} \int^\oplus_{\nu \in \mathfrak{a}_q^*} \pi_{\mu, i \nu}(x) D_{\mu, i\nu}^{-1} dm_\mu(\nu) 
$$
for $ x \in \DF(G_q) $, where 
$$ 
D_{\mu, i \nu} = \pi_{\mu, i\nu}(K_\rho \bowtie 1),
$$
and the measures $ dm_\mu $ on $ \mathfrak{a}_q^* $ are given by 
\begin{align*}
dm_\mu(\nu) &= \frac{1}{|W|} \prod_{\alpha \in {\bf \Delta}^+} |q^{\frac{1}{2}(\alpha, \mu + i \nu)} - q^{-\frac{1}{2}(\alpha, \mu + i \nu)}|^2 \, d\nu. 
\end{align*}
Here $ d\nu $ denotes Lebesgue measure on the torus $ \mathfrak{a}_q^* $, normalized such that its total mass equals $ 1 $.
\end{theorem} 

For $ G_q = \SL_q(2, \mathbb{C}) $ the formulas in Theorem \ref{thmcqgplancherel} reduce to the result obtained by Buffenoir 
and Roche \cite{BRLorentz}. Buffenoir and Roche work with a different normalization of Lebesgue measure and a different parametrization of the irreducible 
representations, so that the formulas given in \cite{BRLorentz} look slightly different. 

We will prove Theorem \ref{thmcqgplancherel} by establishing the Plancherel formula 
$$
\varepsilon_{G_q}(f) = \sum_{\mu \in \weights} \int_{\mathfrak{a}_q^*} \tr(\pi_{\mu, i \nu}(f) D_{\mu, i\nu}^{-2}) dm_\mu(\nu)
$$
for elements of the form $ f = u^\beta_{ij} \otimes \omega^\gamma_{kl} \in \CF^\infty_c(G_q) $ with $ \beta, \gamma \in \weights^+ $. 
Since these elements span $ \CF^\infty_c(G_q) $ linearly this will yield the claim due to Lemma \ref{plancherelformula}. 
For notational convenience let us write 
$$
\tau(f) = \sum_{\mu \in \weights} \int_{\mathfrak{a}_q^*} \tr(\pi_{\mu, i \nu}(f) D_{\mu, i\nu}^{-2}) dm_\mu(\nu)
$$
for the right hand side of the Plancherel formula. We will refer to the functional $ \tau: \CF^\infty_c(G_q) \rightarrow \mathbb{C} $ 
as the \emph{Plancherel functional}. 

\begin{lemma} \label{tauformula}
Let $ f = u^\beta_{ij} \otimes \omega^\gamma_{kl} \in \CF^\infty_c(G_q) $. Then 
\begin{align*}
\tau(f) &= \dim_q(V(\gamma)) 
\sum_{w \in W} \sum_{\varepsilon_r = w \rho - \rho} \sum_m (-1)^{l(w)} 
\phi(u^\gamma_{lr} S^{-1}(u^\beta_{mj}) S^{-1}(u^\gamma_{rk}) u^\beta_{im}) q^{(\varepsilon_r, \varepsilon_m)}, 
\end{align*} 
where $ l(w) $ denotes the length of $ w \in W $. 
\end{lemma} 

\proof From Theorem \ref{thmprincipalirreducible} and Theorem \ref{thmweylgroupactionirreds} we obtain 
$$
\tr(\pi_{y \mu, i y \nu}(f) D_{y \mu, i y \nu}^{-2}) = \tr(\pi_{\mu, i\nu}(f) D_{\mu, i\nu}^{-2}) 
$$
for all $ y \in W $, taking into account that the 
Duflo-Moore operators are given by the action of a certain element of $ \M(\DF(G_q)) $ and so commute with the intertwiners between unitary principal series representations.
Using the Weyl denominator formula 
\begin{align*}
&\prod_{\alpha \in {\bf \Delta}^+} |q^{\frac{1}{2}(\alpha, \mu + i \nu)} - q^{-\frac{1}{2}(\alpha, \mu + i \nu)}|^2 \, 
= \sum_{x, y \in W} (-1)^{l(x) + l(y)} q^{(x \rho + y \rho, \mu)} q^{(x \rho - y \rho, i \nu)}
\end{align*} 
we therefore obtain 
\begin{align*}
&\tau(f) = \frac{1}{|W|} \sum_{\mu \in \weights} \int_{\mathfrak{t}_q} 
\tr(\pi_{\mu, i\nu}(f) D_{\mu, i\nu}^{-2}) 
\prod_{\alpha \in {\bf \Delta}^+} |q^{\frac{1}{2}(\alpha, \mu + i \nu)} - q^{-\frac{1}{2}(\alpha, \mu + i \nu)}|^2 \, d\nu \\
&= \frac{1}{|W|} \sum_{\mu \in \weights} \int_{\mathfrak{t}_q} \sum_{x, y \in W} 
\tr(\pi_{\mu, i\nu}(f) D_{\mu, i\nu}^{-2}) 
(-1)^{l(x) + l(y)} q^{(x \rho + y \rho, \mu)} q^{(x \rho - y \rho, i \nu)} d\nu \\
&= \frac{1}{|W|} \sum_{\mu \in \weights} \int_{\mathfrak{t}_q} \sum_{x, y \in W} 
\tr(\pi_{y \mu, i y \nu}(f) D_{y \mu, i y \nu}^{-2}) 
(-1)^{l(x) + l(y)} q^{(x \rho + y \rho, y \mu)} q^{(x \rho - y \rho, i y \nu)} d\nu \\
&= \frac{1}{|W|} \sum_{\mu \in \weights} \int_{\mathfrak{t}_q} \sum_{x, y \in W} 
\tr(\pi_{\mu, i\nu}(f) D_{\mu, i\nu}^{-2}) 
(-1)^{l(x) + l(y)} q^{(y^{-1} x \rho + \rho, \mu)} q^{(y^{-1}x \rho - \rho, i \nu)} d\nu \\
&= \sum_{\mu \in \weights} \int_{\mathfrak{t}_q} \sum_{w \in W} 
\tr(\pi_{\mu, i\nu}(f) D_{\mu, i\nu}^{-2}) 
(-1)^{l(w)} q^{(w \rho + \rho, \mu)} q^{(w \rho - \rho, i \nu)} d\nu \\
&= \sum_{\mu \in \weights} \int_{\mathfrak{t}_q} \sum_{w \in W} q^{(2 \rho, \mu)} 
\tr(\pi_{\mu, i\nu}(f) D_{\mu, i\nu}^{-2}) 
(-1)^{l(w)} q^{(w \rho - \rho, \mu)} q^{(w \rho - \rho, i \nu)} d\nu.  
\end{align*} 
Inserting the formula 
\begin{align*}
\tr(\pi_{\mu, i \nu}(f) &D_{\mu, i \nu}^{-2}) \\
&=q^{(-2 \rho, \mu)} \dim_q(V(\gamma)) \sum_{\varepsilon_m = \mu} \sum_r  
\phi(u^\gamma_{lr} S^{-1}(u^\beta_{mj}) S^{-1}(u^\gamma_{rk}) u^\beta_{im}) q^{-(i \nu, \varepsilon_r)} 
\end{align*}
for the twisted character 
from Proposition \ref{characterformula} we arrive at the claimed expression
for the Plancherel functional. \qed 

Recall next the construction of the BGG complex for quantized universal enveloping algebras from \cite{HeckenbergerKolbBGG}. 
As mentioned in the introduction, the quantum BGG complex in its original form is a resolution of a finite dimensional representation of $U_q(\mathfrak{g})$ by direct sums of Verma modules, generalizing a resolution for classical enveloping algebras first discovered in \cite{BGG_differentialoperators}.  
Specifically, for a dominant integral 
weight $ \nu \in \weights^+ $ and $ k \geq 0 $, set 
$$
C_k = \bigoplus_{\substack{w \in W \\ l(w) = k}} M(w . \nu),  
$$
where $ M(\eta) $ denotes the Verma module of $ U_q(\mathfrak{g}) $ with highest weight $ \eta \in \weights $ and 
$$
w . \nu = w(\nu + \rho) - \rho 
$$ 
is the shifted Weyl group action. 
Using inclusions of Verma modules, one constructs boundary operators $ d: C_k \rightarrow C_{k - 1} $ such that $ d^2 = 0 $ in the same way as for the 
original BGG complex. Let us also denote by $ \varepsilon_\nu: C_0 = M(\nu) \rightarrow V(\nu) $ the canonical projection, where $ V(\nu) $ is the unique irreducible 
quotient of $ M(\nu) $. 

The following theorem is a special case of the results obtained by Heckenberger and Kolb in \cite{HeckenbergerKolbBGG}. 

\begin{theorem} \label{HKBGG}
The chain complex 
$$
\xymatrix{
0 \ar@{->}[r] & C_n \ar@{->}[r]^d & C_{n - 1} \ar@{->}[r]^d & \cdots \ar@{->}[r]^d & C_0 \ar@{->}[r]^{\varepsilon_\nu} & V(\nu) \ar@{->}[r] & 0
}
$$
is exact.  
\end{theorem} 

We shall use the category equivalence between category $ \O $ and Harish-Chandra modules \cite{Josephbook}, \cite{VYcqg} to transport the BGG complex 
from Theorem \ref{HKBGG} for $ \nu = 0 $ into a complex of principal series modules. 

Following the notation and terminology in \cite[Section 6.5]{VYcqg}, let us fix $ l = 0 $ and denote by $ \O_l $ the full subcategory of category $ \O $ consisting 
of modules whose weights all belong to $ \weights \subset \mathfrak{h}^*_q $. Also, let $ \HC_l $ be the full subcategory of 
Harish-Chandra bimodules for which the annihilator of the right action of $ ZU_q(\mathfrak{g}) $ contains the kernel of the central character associated 
with $ l $. Since the weight $ l = 0 $ is dominant and regular the 
functor $ \F_l: \O_l \rightarrow \HC_l $ defined by 
$$
\F_l(M) = F \Hom(M(l), M)
$$
is an equivalence of categories, see Section 6.5 in \cite{VYcqg}. 
Using duality in category $ \O $, we conclude that setting $ D^k = \F_l(C_k^\vee) $ and $ \partial = \F_l(d^\vee) $ yields an exact complex of Harish-Chandra 
modules.  

More explicitly, setting 
$$
\mu = l - r, \qquad \lambda + 2 \rho = -l - r
$$
where $ l = 0 $ and $ r = w . 0 = w \rho - \rho $ we get an isomorphism $ F \Hom(M(l), M(r)^\vee) \cong \Gamma(\E_{\mu, \lambda}) $ of $ \DF(G_q) $-modules, 
so that
$$ 
\F_l(M(w . 0)^\vee) \cong \Gamma(\E_{-w . 0, -w . 0 - 2 \rho}). 
$$
Remark that $ \F_l(M(0)^\vee) = \Gamma(\E_{0, -2\rho}) $ contains the trivial $ \DF(G_q) $-module $ \mathbb{C} \cong \F_l(V(0)^\vee) $ as a submodule.   
We thus arrive at the exact chain complex 
$$
\xymatrix{
0 \ar@{->}[r] & \mathbb{C} \ar@{->}[r]^\iota & D^0 \ar@{->}[r]^\partial & D^1 \ar@{->}[r]^\partial & \cdots \ar@{->}[r]^\partial & D^n \ar@{->}[r] & 0
}
$$
of $ \DF(G_q) $-modules, where 
$$
D^j = \bigoplus_{\substack{w \in W \\ l(w) = j}} \Gamma(\E_{-w . 0, -w . 0 - 2\rho}). 
$$
We will refer to $ D^\bullet $ as the \emph{geometric BGG complex}. Let us point out that $ D^\bullet $ is naturally a cochain complex, so that we are 
using cohomological indexing. 

Now assume that $ p = p \bowtie 1 \in \DF(K_q) \bowtie \CF^\infty(K_q) = \DF(G_q) $ is a finite projection, 
by which we mean that $ p $ is supported on finitely many matrix blocks $ \KH(V(\nu)) $ inside $ \DF(K_q) $. 
Then the action of $ p $ on the spaces $ D^j $ determines direct summands $ p \cdot D^j \subset D^j $ which assemble into a 
subcomplex $ (p \cdot D)^\bullet $ of the geometrical BGG complex. Observe that $ (p \cdot D)^\bullet $ is in fact a finite dimensional exact 
complex of $ p \DF(G_q) p $-modules since all isotpyical components of principal series modules are finite dimensional. 

Recall the following basic fact from homological algebra. Assume that $ C^\bullet $ is a finite dimensional complex of vector spaces, so that 
all the spaces $ C^n $ are finite dimensional and $ C^\bullet $ is supported in finitely many degrees. Moreover 
let $ f: C^\bullet \rightarrow C^\bullet $ be a chain map, with induced maps 
$ H^*(f): H^*(C) \rightarrow H^*(C) $ on cohomology. Since the cohomology groups of $ C^\bullet $ are finite dimensional as well we can form 
$ \tr(H^k(f)) $, and we have the Hopf trace formula 
$$
\sum_{k \in \mathbb{Z}} (-1)^k \tr_{H^k(C)}(H^k(f)) = \sum_{k \in \mathbb{Z}} (-1)^k \tr_{C^k}(f^k), 
$$ 
where $ \tr_V $ denotes the natural trace on a vector space $ V $. 

Combining the above considerations we arrive at the following key lemma. 

\begin{lemma} \label{BGGhopftraceformula}
Let $ x \in \DF(G_q) $. Then 
$$
\sum_{w \in W} (-1)^{l(w)} \tr(\pi_{-w . 0, -w . 0 - 2 \rho}(x)) = \hat{\varepsilon}_{G_q}(x),  
$$
where $ \hat{\varepsilon}_{G_q} $ is the counit of $ \DF(G_q) $. 
\end{lemma} 

\proof Let us first point out that the operators $ \pi_{-w . 0, -w . 0 - 2 \rho}(x) $ are finite rank, so that the left hand side 
of the above formula is well-defined. 

Since the boundary operators in the geometric BGG complex are $ \DF(G_q) $-linear the endomorphism of $ D^\bullet $ induced 
by $ x $ is a chain map. Note in addition that we can find a finite central projection $ p = p \bowtie 1 \in \DF(G_q) = \DF(K_q) \bowtie \CF^\infty(K_q) $ 
such that $ x = p x p $. Hence the operator trace of $ \pi_{-w . 0, -w . 0 - 2 \rho}(x) $ equals the trace of $ \pi_{-w . 0, -w . 0 - 2 \rho}(pxp) $ viewed 
as endomorphism of the finite dimensional vector space $ p \cdot \Gamma(\E_{-w . 0, -w . 0 - 2 \rho}) $. 
Now the Hopf trace formula applied to the complex $ (p \cdot D)^\bullet $ yields the claim, 
using that the action of $ x $ on the trivial $ \DF(G_q) $-module $ \mathbb{C} \subset D^0 $ is given by the counit $ \hat{\varepsilon}_G $. \qed 

Let us now go back to the problem of computing the Plancherel functional $ \tau(f) $ for an element of the 
form $ f = u^\beta_{ij} \otimes \omega^\gamma_{kl} \in \CF^\infty_c(G_q) $. The following discussion completes the proof of 
Theorem \ref{thmcqgplancherel}. 

\begin{theorem} \label{plancherelmain}
Let $ f = u^\beta_{ij} \otimes \omega^\gamma_{kl} \in \CF^\infty_c(G_q) $. Then $ \tau(f) = \varepsilon_{G_q}(f) $. 
\end{theorem} 

\proof According to Lemma \ref{tauformula} we can write
$$
\tau(f) = \sum_{w\in W} (-1)^{l(w)} \tau_w(f)
$$
where
\begin{align*}
&\tau_w(f) = \dim_q(V(\gamma)) \sum_{\varepsilon_r = w.0} \sum_m
\phi(u^\gamma_{lr} S^{-1}(u^\beta_{mj}) S^{-1}(u^\gamma_{rk}) u^\beta_{im}) q^{(\varepsilon_r,\varepsilon_m)} \\
&= \dim_q(V(\gamma)) \sum_{\varepsilon_r = w.0} \sum_{m,n}
\phi(u^\gamma_{lr} S^{-1}(u^\beta_{nj}) S^{-1}(u^\gamma_{rk}) S(S^{-1}(u^\beta_{im})) (K_{-w.0}, S^{-1}(u^\beta_{mn}))) \\
&= \dim_q(V(\gamma)) \sum_{\varepsilon_r = w.0}
\phi(u^\gamma_{lr} \pi_{-w.0,-w.0-2\rho}(1 \bowtie S^{-1}(u^\beta_{ij}))(S^{-1}(u^\gamma_{rk}))), 
\end{align*}
using the definition of the Yetter-Drinfeld action of $ \CF^\infty(K_q) \subset \M(\DF(G_q)) $ on the principal series module $ \Gamma(\E_{-w.0,-w.0 - 2\rho}) $ 
in the last step.

Note that the vectors $ e^\nu_{ab} = S^{-1}(u^\nu_{ba}) $ with $ \nu \in \weights^+, \varepsilon_b = w . 0 $ and $ a $ arbitrary form a basis 
of $ \Gamma(\E_{-w.0,-w.0 - 2\rho}) $. If we consider the linear functionals $ e_\nu^{ab} $ on $ \Gamma(\E_{-w.0,-w.0 - 2\rho}) $ defined by 
$$ 
e_\nu^{ab}(\xi) = q^{(2 \rho, \varepsilon_a)} \dim_q(V(\nu)) \phi(u^\nu_{ab} \xi), 
$$ 
then we obtain 
\begin{align*}
e_\eta^{cd}(e^\nu_{ab}) &= q^{(2 \rho, \varepsilon_c)} \dim_q(V(\nu)) \phi(u^\eta_{cd} S^{-1}(u^\nu_{ba})) \\
&= q^{(2 \rho, \varepsilon_c)} \dim_q(V(\nu)) \phi(u^\nu_{ba} S(u^\eta_{cd})) \\
&= \delta_{\eta \nu} \delta_{ac} \delta_{bd} 
\end{align*}
for any $ \eta, \nu \in \weights^+ $ and $ a,b,c,d $ according to the Schur orthogonality relations. 
It follows that the vectors $ e_\eta^{cd} $ are the dual basis vectors to the vectors $ e^\nu_{ab} $. 

Observe next that the action of $ \hat{S}^{-2}(\omega^\gamma_{kl}) \bowtie 1 \in \DF(G_q)$ on $ \Gamma(\E_{-w.0,-w.0 - 2\rho}) $ is given by 
\begin{align*}
\pi_{-w.0,-w.0-2\rho}(\hat{S}^{-2}(\omega^\gamma_{kl}) \bowtie 1) (e^\eta_{sr}) 
&= \pi_{-w.0,-w.0-2\rho}(\hat{S}^{-2}(\omega^\gamma_{kl}) \bowtie 1) (S^{-1}(u^\eta_{rs})) \\
&= \sum_t (\hat{S}^{-1}(\omega^\gamma_{kl}), S^{-1}(u^\eta_{ts})) S^{-1}(u^\eta_{rt}) \\
&= \delta_{\gamma \eta} \delta_{sl} S^{-1}(u^\eta_{rk}) = \delta_{\gamma \eta} \delta_{sl} e^\eta_{kr}.
\end{align*}

We shall now assemble these considerations. More precisely, let us consider the element $ x \in \DF(G_q) $ defined by 
$$ 
x = q^{(-2\rho, \varepsilon_l)} \hat{S}^{-2}(\omega^\gamma_{kl}) \bowtie S^{-1}(u^\beta_{ij}) 
= q^{(-2\rho, \varepsilon_l)} (\hat{S}^{-2}(\omega^\gamma_{kl}) \bowtie 1)(1 \bowtie S^{-1}(u^\beta_{ij})). 
$$ 
Then by combining the above formulas we compute 
\begin{align*}
\tr(&\pi_{-w.0,-w.0-2\rho}(x)) = \sum_{\eta \in \weights^+} \sum_{\varepsilon_r = w.0} \sum_s e_\eta^{sr}(\pi_{-w.0,-w.0-2\rho}(x)(e^\eta_{sr})) \\
&= q^{(-2 \rho, \varepsilon_l)} \sum_{\varepsilon_r = w.0} e_\gamma^{lr}(\pi_{-w.0,-w.0-2\rho}(1 \bowtie S^{-1}(u^\beta_{ij}))(e^\gamma_{kr})) \\
&= \dim_q(V(\gamma)) \sum_{\varepsilon_r = w.0} \phi(u^\gamma_{lr} \pi_{-w.0,-w.0-2\rho}(1 \bowtie S^{-1}(u^\beta_{ij}))(S^{-1}(u^\gamma_{rk}))) \\
&= \tau_w(f). 
\end{align*}
Applying the Hopf trace formula from Lemma \ref{BGGhopftraceformula}, we conclude that $ \tau(f) $ is equal to $ \hat{\varepsilon}_{G_q}(x) $. 
In other words, we obtain
$$
\tau(f) = \hat{\varepsilon}_{G_q}(x) = \delta_{ij} \delta_{\gamma 0} = \varepsilon_{G_q}(f) 
$$
as desired. \qed 

Let us remark that the Plancherel measure for $ G_q $ resembles its counterpart for the classical group $ G $. More precisely, up to a normalization 
the classical measure is given by 
$$
\prod_{\alpha \in {\bf \Delta}^+} |(\alpha, \mu + i \nu)|^2 d\nu
= \prod_{\alpha \in {\bf \Delta}^+} (\alpha, \mu + i \nu) (\alpha, \mu - i \nu) d\nu
$$
in the component of the parameter space $ \weights \times \mathfrak{a}^* $ corresponding to $ \mu \in \weights $, compare Section 5 
in \cite{HarishChandraplancherelcomplex}. 
By comparison, the measure in Theorem \ref{thmcqgplancherel} reads 
\begin{align*}
dm_\mu(\nu) &= \frac{1}{|W|} \prod_{\alpha \in {\bf \Delta}^+} |q^{\frac{1}{2}(\alpha, \mu + i \nu)} - q^{-\frac{1}{2}(\alpha, \mu + i \nu)}|^2 \, d\nu \\
&= \frac{1}{|W|} \prod_{\alpha \in {\bf \Delta}^+} (q - q^{-1})^2 
[\tfrac{1}{2}(\alpha, \mu + i \nu)]_q[\tfrac{1}{2}(\alpha, \mu - i \nu)]_q \, d\nu. 
\end{align*}
Expanding $ q = e^h $ in powers of $ h $ this can be rewritten as  
\begin{align*}
dm_\mu(\nu) 
&= \frac{1}{|W|} \prod_{\alpha \in {\bf \Delta}^+} (q^{(\alpha, \mu)} + q^{-(\alpha, \mu)} - q^{(\alpha, i\nu)} - q^{-(\alpha, i\nu)}) d\nu \\
&= \frac{1}{|W|} \prod_{\alpha \in {\bf \Delta}^+} h^2 (\alpha, \mu + i \nu) (\alpha, \mu - i \nu) d\nu + \text{ higher order terms}. 
\end{align*}
That is, up to a scalar, the first nonzero coefficient in the expansion agrees with the formula for the classical measure.

\section{Further remarks on the main result} \label{secproofslq2} 

In this section we include a few supplementary remarks on computational aspects of the formulas obtained in Theorem \ref{plancherelmain}. 

Let $ \nu \in \weights^+ $ and let $ e^\nu_1, \dots, e^\nu_n $ be an orthonormal weight basis of $ V(\nu) $. Moreover denote by $ e_\nu^1, \dots, e_\nu^n $ 
the dual basis of the contragredient representation $ V(\nu)^* = \Hom(V(\nu), \mathbb{C}) $. 
Inspecting Lemma \ref{tauformula} and the definition of the Haar functional, the computation of the Plancherel functional $ \tau(f) $ for all $ f $ of the 
form $ f = u^\beta_{ij} \otimes \omega^\gamma_{kl} $ with fixed $ \beta, \gamma \in \weights^+ $ can be rephrased in terms of the tensor 
$$
\tau_{\beta \gamma} = \sum_{w \in W} \sum_{\varepsilon_r = w \rho - \rho} \sum_m (-1)^{l(w)} q^{(\varepsilon_r, \varepsilon_m - 2 \rho)} 
P(e^\gamma_r \otimes e^\beta_m \otimes e_\gamma^r \otimes e_\beta^m), 
$$
where $ P $ denotes the orthogonal projection onto the trivial isotypical component of the tensor 
product $ V(\gamma) \otimes V(\beta) \otimes V(\gamma)^* \otimes V(\beta)^* $. 
More precisely, the nontrivial part of Theorem \ref{plancherelmain}, namely the vanishing of all $ \tau(f) $ for all $ f = u^\beta_{ij} \otimes \omega^\gamma_{kl} $ 
with $ \gamma \neq 0 $, is equivalent to the following assertion. 

\begin{theorem} \label{plancherelmaininvariants}
For all $ \beta, \gamma \in \weights^+ $ with $ \gamma \neq 0 $ we have $ \tau_{\beta \gamma} = 0 $. 
\end{theorem} 

Except in a few special cases, it seems a forbidding task to compute any of the summands 
appearing in the tensor $ \tau_{\beta \gamma} $ explicitly. 

However, let us restrict attention to the case of the quantum Lorentz group $ G_q = \SL_q(2, \mathbb{C}) $ and explain how to verify 
Theorem \ref{plancherelmaininvariants} by elementary calculations in this case nonetheless. This yields a shorter proof of the Plancherel formula 
than the original one by Buffenoir and Roche \cite{BRLorentz}, and does not invoke any homological algebra arguments. 

We identify the set of weights $ \weights $ of $ K_q = \SU_q(2) $ with $ \frac{1}{2} \mathbb{Z} $.
Moreover we shall work with the orthonormal basis $ e^\nu_j $ for $ j \in \{-\nu, -\nu + 1, \dots, \nu \} $ of 
the irreducible representation $ V(\nu) $ of heighest weight $ \nu \in \frac{1}{2} \mathbb{N}_0 $ as in \cite[Section 6.8.3]{VYcqg}. 
Explicitly, we have 
\begin{align*}
E \cdot e^\nu_j &= q^j [\nu - j]_q^{\frac{1}{2}} [\nu + j + 1]_q^{\frac{1}{2}} e^\nu_{j + 1}, \\
F \cdot e^\nu_j &= q^{-(j - 1)} [\nu + j]_q^{\frac{1}{2}} [\nu - j + 1]_q^{\frac{1}{2}} e^\nu_{j - 1},
\end{align*}
where we interpret $ e^\nu_j = 0 $ if $ |j| > \nu $, and we abbreviate $ E = E_1, F = F_1 $. 
For the dual basis vectors $ e_\nu^j \in V(\nu)^* $ in the contragredient representation we obtain 
\begin{align*}
E \cdot e_\nu^j &= - q^{-j + 1} [\nu - j + 1]_q^{\frac{1}{2}} [\nu + j]_q^{\frac{1}{2}} e_\nu^{j - 1}, \\
F \cdot e_\nu^j &= -q^{j} [\nu + j + 1]_q^{\frac{1}{2}} [\nu - j]_q^{\frac{1}{2}} e_\nu^{j + 1}. 
\end{align*}
Using these formulas we shall verify the following relation, where we write again $ P $ for the projection onto the 
trivial isotypical component. 

\begin{lemma} \label{invlemma} 
Let $ \beta, \gamma \in \frac{1}{2} \mathbb{N}_0 $. If $ \gamma > 0 $ then 
\begin{align*}
\sum_m q^{-2 m(r + 1) + 2} P(e^\gamma_{r - 1} \otimes e^\beta_m \otimes e_\gamma^{r - 1} \otimes e_\beta^m) 
= \sum_m q^{-2 mr}  P(e^\gamma_r \otimes e^\beta_m \otimes e_\gamma^r \otimes e_\beta^m) 
\end{align*}
for all $ r \in \{-\gamma + 1, -\gamma + 2, \dots, \gamma\} $. 
\end{lemma} 

\proof Let us consider the relation 
$$ 
0 = \sum_{m \in \weights} q^{-2(m - 1)r} P(E \cdot (e^\gamma_{r - 1} \otimes e^\beta_m \otimes e_\gamma^r \otimes e_\beta^m)), 
$$ 
obtained from the fact that $ E $ acts by zero on the trivial representation. We calculate
\begin{align*}
q^{-2(m - 1)r} &(E \cdot e^\gamma_{r - 1}) \otimes (K^2 \cdot e^\beta_m) \otimes (K^2 \cdot e_\gamma^r) \otimes (K^2 \cdot e_\beta^m) \\
&= q^{-2(m - 1)r} q^{-2r} q^{r - 1} [\gamma - r + 1]_q^{\frac{1}{2}}[\gamma + r]_q^{\frac{1}{2}} e^\gamma_r \otimes e^\beta_m \otimes e_\gamma^r \otimes e_\beta^m \\
&= q^{-2mr} q^{r - 1} [\gamma - r + 1]_q^{\frac{1}{2}}[\gamma + r]_q^{\frac{1}{2}} e^\gamma_r  \otimes e^\beta_m \otimes e_\gamma^r \otimes e_\beta^m, 
\end{align*} 
using the notation $ K^2 = K_1 $. Similarly, 
\begin{align*}
&q^{-2(m - 1)r} (1 \cdot e^\gamma_{r - 1}) \otimes (1 \cdot e^\beta_m) \otimes (E \cdot e_\gamma^r) \otimes (K^2 \cdot e_\beta^m) \\
&= - q^{-2(m - 1)r} q^{-2 m} q^{-r + 1} [\gamma - r + 1]_q^{\frac{1}{2}} [\gamma + r]_q^{\frac{1}{2}}
e^\gamma_{r - 1} \otimes e^\beta_m \otimes e_\gamma^{r - 1} \otimes e_\beta^m \\
&= - q^{-2 m(r + 1)} q^{r + 1} [\gamma - r + 1]_q^{\frac{1}{2}} [\gamma + r]_q^{\frac{1}{2}}
e^\gamma_{r - 1} \otimes e^\beta_m \otimes e_\gamma^{r - 1} \otimes e_\beta^m. 
\end{align*}
Moreover we have 
\begin{align*}
&\sum_m q^{-2(m - 1)r}
(1 \cdot e^\gamma_{r - 1}) \otimes (E \cdot e^\beta_m) \otimes (K^2 \cdot e_\gamma^r) \otimes (K^2 \cdot e_\beta^m) \\
&\quad + q^{-2(m - 1)r} (1 \cdot e^\gamma_{r - 1}) \otimes (1 \cdot e^\beta_m) \otimes (1 \cdot e_\gamma^r) \otimes (E \cdot e_\beta^m) \\
&= \sum_m q^{-2(m - 1)r} q^{-2r} q^{-2 m} q^m [\beta - m]_q^{\frac{1}{2}} [\beta + m + 1]_q^{\frac{1}{2}} 
e^\gamma_{r - 1} \otimes e^\beta_{m + 1} \otimes e_\gamma^r \otimes e_\beta^m \\
&\quad - q^{-2(m - 1)r} q^{-m + 1}[\beta - m + 1]_q^{\frac{1}{2}} [\beta + m]_q^{\frac{1}{2}} 
e^\gamma_{r - 1} \otimes e^\beta_m \otimes e_\gamma^r \otimes e_\beta^{m - 1} \\
&= \sum_m q^{-2(m - 1)r} q^{-2r} q^{-m} [\beta - m]_q^{\frac{1}{2}} [\beta + m + 1]_q^{\frac{1}{2}} 
e^\gamma_{r - 1} \otimes e^\beta_{m + 1} \otimes e_\gamma^r \otimes e_\beta^m \\
&\quad - q^{-2mr} q^{-m} [\beta - m]_q^{\frac{1}{2}} [\beta + m + 1]_q^{\frac{1}{2}} 
e^\gamma_{r - 1} \otimes e^\beta_{m + 1} \otimes e_\gamma^r \otimes e_\beta^m \\
&= 0. 
\end{align*} 
Combining these relations yields the claim. \qed 

We may now compute $ \tau_{\beta \gamma} $ for all $ \beta, \gamma \in \frac{1}{2} \mathbb{N}_0 $ with $ \gamma > 0 $. More precisely, according to 
Lemma \ref{invlemma} we obtain 
\begin{align*}
\tau_{\beta \gamma} = \sum_m P(e^\gamma_0 \otimes e^\beta_m \otimes e_\gamma^0 \otimes e_\beta^m) 
- \sum_m q^{-2m + 2} P(e^\gamma_{-1} \otimes e^\beta_m \otimes e_\gamma^{-1} \otimes e_\beta^m) = 0. 
\end{align*}
This proves Theorem \ref{plancherelmaininvariants} in rank $ 1 $, and 
therefore also the Plancherel formula for the quantum group $ G_q = \SL_q(2, \mathbb{C}) $.

\section{The reduced $ C^* $-algebras of complex quantum groups} \label{secreducedcstar}

In this section we use the Plancherel Theorem \ref{thmcqgplancherel} to describe the structure of the reduced group $ C^* $-algebras of complex 
quantum groups, in analogy with the classical case. 

Let $ G_q $ be a complex semisimple quantum group. 
Moreover let $ \H = (\H_{\mu, \lambda})_{\mu, \lambda} $ be the locally constant Hilbert space bundle of unitary principal series representations of $ G_q $ 
over $ \weights \times \mathfrak{t}_q^* $. By slight abuse of notation we will also write $ \H $ for the corresponding 
Hilbert $ C_0(\weights \times \mathfrak{t}_q^*) $-module. 
Inspecting the explicit formulas for the action of $ \DF(G_q) $ on unitary principal series representations we obtain a 
nondegenerate $ * $-homomorphism $ \pi: \DF(G_q) \rightarrow C_0(\weights \times \mathfrak{t}_q^*, \KH(\H)) $ by 
setting $ \pi(x)(\mu, \lambda) = \pi_{\mu, \lambda}(x) $. 
Here $ C_0(\weights \times \mathfrak{t}_q^*, \KH(\H)) $ denotes the $ C_0(\weights \times \mathfrak{t}_q^*) $-algebra of compact operators on 
the Hilbert $ C_0(\weights \times \mathfrak{t}_q^*) $-module $ \H $. 
By the definition of the maximal group $ C^* $-algebra $ C^*_\max(G_q) $, the map $ \pi $ extends uniquely to a 
nondegenerate $ * $-homomorphism $ C^*_\max(G_q) \rightarrow C_0(\weights \times \mathfrak{t}_q^*, \KH(\H)) $, which will again be denoted by $ \pi $, and 
which we will refer to as the \emph{canonical $ * $-homomorphism} below. 

We obtain an action of $ W $ on $ C_0(\weights \times \mathfrak{t}_q^*, \KH(\H)) $ by 
$$ 
(w \cdot f)(\mu, \lambda) = U(w)_{\mu, \lambda} f(w^{-1} \mu, w^{-1} \lambda) U(w)_{\mu, \lambda}^*, 
$$
where $ U(w)_{\mu, \lambda}: \H_{w^{-1} \mu, w^{-1} \lambda} \rightarrow \H_{\mu, \lambda} $ is a unitary intertwiner as in Theorem \ref{thmweylgroupactionirreds}. 
Using the explicit construction of intertwiners associated to simple reflections in Section 6.8 of \cite{VYcqg} one checks that 
this is indeed a well-defined action, noting that conjugation by $ U(w)_{\mu, \lambda} $ eliminates the scalar ambiguity in the definition 
of the intertwiners. 

With this notation and terminology in place the structure of $ C^*_\red(G_q) $ can be described as follows. 

\begin{theorem} \label{reducedcstarstructure}
Let $ G_q $ be a complex semisimple quantum group, and let $ \H = (\H_{\mu, \lambda})_{\mu, \lambda} $ be the 
Hilbert space bundle of unitary principal series representations of $ G_q $ over $ \weights \times \mathfrak{t}_q^* $. Then one obtains an isomorphism 
$$
C^*_\red(G_q) \cong C_0(\weights \times \mathfrak{t}_q^*, \KH(\H))^W 
$$
induced by the canonical $ * $-homomorphism $ \pi: C^*_\max(G_q) \rightarrow C_0(\weights \times \mathfrak{t}_q^*, \KH(\H)) $. 
\end{theorem}

\proof According to the Plancherel Theorem \ref{thmcqgplancherel}, all unitary principal series representations of $ G_q $ factorize over the 
reduced group $ C^* $-algebra, and the resulting $ * $-homomorphism $ \pi: C^*_\red(G_q) \rightarrow C_0(\weights \times \mathfrak{t}_q^*, \KH(\H)) $ 
is injective. The image $ \im(\pi) $ of the map $ \pi $ is contained in $ C_0(\weights \times \mathfrak{t}_q^*, \KH(\H))^W $ by construction. 

It remains to show that $ \im(\pi) $ is in fact equal to $ C_0(\weights \times \mathfrak{t}_q^*, \KH(\H))^W $. 
Note that the irreducible representations of $ A = C_0(\weights \times \mathfrak{t}_q^*, \KH(\H))^W $ 
are given by point evaluations on $ \weights \times \mathfrak{t}_q^* $, 
and these remain irreducible when restricted to the image of $ \pi $. Moreover, two irreducible representations of $ A $ are inequivalent iff they 
correspond to parameters in different orbits of the Weyl group action on $ \weights \times \mathfrak{t}_q^* $. 
According to Theorem \ref{thmweylgroupactionirreds} the same condition distinguishes unitary principal series representations. 
Hence Dixmier's version of the Stone-Weierstrass Theorem, see Section 11.1 in \cite{Dixmiercstar}, yields the claim. \qed 

Theorem \ref{reducedcstarstructure} shows in particular that the trivial representation of $ G_q $ does not factorize through $ C^*_\red(G_q) $. 
In other words, the full and reduced group $ C^* $-algebras of $ G_q $ are not isomorphic, which means that $ G_q $ is not amenable \cite{BedosTuset}. 
More interestingly, Arano has shown \cite{Aranospherical}, \cite{Aranocomparison} that higher rank complex quantum groups do in fact have property (T). 

Finally, let us point out that Theorem \ref{reducedcstarstructure} illustrates nicely the deformation aspect of the theory of complex semisimple quantum groups, 
a feature which is not apparent from the Drinfeld double construction. Indeed, by setting formally $ h = 0 $ and $ \mathfrak{t}^*_1 = \mathfrak{t}^* $ in 
Theorem \ref{reducedcstarstructure} we reobtain the well-known description of the reduced group $ C^* $-algebra of the classical complex semisimple Lie group $ G $. 
Thus the limit $ q \rightarrow 1 $ corresponds to the opening of the torus $ \mathfrak{t}_q^* = \mathfrak{t}^*/i \hbar^{-1} \roots^\vee $ to $ \mathfrak{t}^* $
as $ \hbar \rightarrow 0 $.

\bibliographystyle{plain}

\bibliography{cvoigt}

\end{document}